\documentclass[11pt]{article}
\renewcommand{\baselinestretch}{1}
\begin{document}

\title{{\bf ASYMPTOTIC BEHAVIOR OF THE NUMBER OF LOST MESSAGES}}
\author{{VYACHESLAV M. ABRAMOV}\thanks{Department of Mathematics,
The Faculty of Exact Sciences, Tel Aviv University, Tel Aviv
69978, Israel and College of Judea and Samaria, Ariel 44837,
Israel. \newline
 E-mail: vyachesl@inter.net.il}} \maketitle

\begin{abstract}
\indent

The goal of the paper is to study asymptotic behavior of the
number of lost messages. Long messages are assumed to be divided
into a random number of packets which are transmitted
independently of one another. An error in transmission of a packet
results in the loss of the entire message. Messages arrive to the
$M/GI/1$ finite buffer model and can be lost in two cases as
either at least one of its packets is corrupted or the buffer is
overflowed. With the parameters of the system typical for models
of information transmission in real networks we obtain theorems on
asymptotic behavior of the number of lost messages. We also study
how the loss probability changes if redundant packets are added.
Our asymptotic analysis approach is based on Tauberian
theorems with remainder.\\



\noindent

\end{abstract}

\section{Introduction}
\renewcommand{\theequation}{1.\arabic{equation}}
\setcounter{section}{1}
\indent

{\bf 1.1. Review of the literature and general description of the
system.} Long messages in Internet protocols, that have to be
transmitted, are divided into small packets. Upon transmission
each packet is transformed by providing additional information
related to a given message. Because of the bit errors in
transmission of the packet, the message can be lost. The loss
probability of a message plays a significant role in the
evaluation of network performance and design of network topology.

There are a number of papers where the loss probability of a
message has been studied. Cidon {\it et al} [11] derived
recurrence relations for the loss probabilities of packets in a
message giving the numerical results for the $M/M/1/n$ buffer
model. The complexity of recurrence calculations of that paper are
$O(nm^2)$ where $m$ is the size of a message, and $n$ is the
buffer capacity. Considering the same model Gurewitz {\it et al}
[13] obtained another representation for the loss probability by
using the ballot theorem (e.g. Tak\'acs [16]). In framework of the
same model Altman and Jean-Marie [7] give a comprehensive analysis
for the multi-dimensional generating function of the loss
probabilities based on the recurrence relations of the paper of
Cidon {\it et al} [11] and analyze the effect of adding redundant
packets. Studying a slightly more general model with several
sources Ait-Hellal {\it et al} [6] obtained some asymptotic
results and studied the effect of adding redundancy to the loss
probability. The aforementioned papers [6], [7], [11], [13] all
discuss the problem of complexity of calculations as well as the
required memory to store intermediate variables.

In real communication networks the capacity is large. Therefore,
asymptotic analysis of the number of lost messages is necessary.
The present paper provides asymptotic analysis with sequential
application to redundancy of the following model. Assume that each
message is divided into a random number of packets each of which
is forwarded to the buffer. For the $i$th message denote its
random number of its packets by $\nu_i$. We assume that the
sequence $\nu_i\ge 1$, $i\ge 1$, consists of independent
identically distributed integer random variables. The interarrival
times between messages have an exponential distribution with
parameter $\lambda$. The buffer can contain only $N$ packets, that
is, if immediately before the arrival of message of $l$ packets
there are $L$ packets in buffer then the message is accepted only
if $L+l\le N$, otherwise the message of $l$ packets is lost. The
loss of a message can also occur if at least one packet in a
message is corrupted. In this case we assume that if there is
enough space then the message does occupy the buffer, but it is
hidden and therefore lost. The probability that at least one
packet in a message is corrupted is denoted by $p$.

In general loss communication networks a transmission time
typically depends on the number of packets in message. To be
realistic we must study a general queueing system with service
time depending on batch size. The analysis of such system is a
hard problem. On the other hand, the model with a fixed number of
packets in message, leading to the standard $M/GI/1/n$ queueing
system, is not realistic. Therefore, in the following we assume
additionally that the random variables $\nu_i$ have fixed upper
and lower bounds $\nu^{upper}$ and $\nu^{lower}$, i.e. ${\bf P}\{
\nu^{lower}\le\nu_i\le\nu^{upper}\}=1$. This assumption can be
considered as a compromise between these two cases.
It has a real application in some communication technologies,
especially in optical local networks, where a number of small
messages following the same direction are combined to one message
(bus)\footnote{For example, one of such technologies was developed
in Orika Optical Networks Limited, where the author worked during
2000-2001 years}. Outgoing from the local network, the bus
continues on its way being processed by the Internet protocols. When
the difference between $\nu^{upper}$ and $\nu^{lower}$ for the
messages is not large, then assumption that a
transmission time is independent of the message size seems
appropriate.

\medskip

{\bf 1.2. Formulation of the model in terms of the queueing
theory.} In terms of the queueing theory the model can be
described as follows. We assume that messages arrive to the finite
buffer $M/GI/1$ queue with random number of waiting places
$\zeta$. The input rate is equal to $\lambda$, the service time
distribution is $B(x)$ with the expectation $b$. By a queueing
system with random number of waiting places we mean the following.
We denote
$$
\zeta = \sup\Big\{m: \sum_{i=1}^m\nu_i\le N\Big\},
$$
and according to the assumption ${\bf P}\{
\nu^{lower}\le\nu_i\le\nu^{upper}\}=1$, there are two fixed values
$\zeta^{upper}$ and $\zeta^{lower}$ depending on $N$, and ${\bf
P}\{ \zeta^{lower}\le\zeta\le\zeta^{upper}\}=1$.

Let $\zeta_1$, $\zeta_2,...$, be a strictly stationary and ergodic
sequence of random variables, ${\bf P}\{\zeta_i=j\}={\bf
P}\{\zeta=j\}$,~$\zeta^{lower}\le j\le\zeta^{upper}$. If
$\xi_i$ is the number of messages in the queue immediately before
arrival of the $i$th message, then the message is lost if
$\xi_i>\zeta_i$. Otherwise it joins the queue. We assume that
$\xi_1=0$.

The existence of the stationary queue-length distribution, i.e.
$$
{\bf P}\{{\mathaccent "016 q}=j\}=\lim_{i\to\infty}{\bf
P}\{\xi_i=j|\xi_1<\infty\},~~~~~j=0,1,...,\zeta^{upper}, \leqno
(1.1)
$$
is shown in the following. The special case when ${\bf
P}\{\nu_i=l\}=1$ leads to the standard $M/GI/1/n$ queueing system,
where $n=\lfloor N/l \rfloor$ is the integer part of $N/l$.

It is also assumed that each message is marked with probability
$p$. We study the asymptotic behavior of the loss probability
under assumptions that ${\bf E}\zeta$ increases to infinity and
$p$ vanishes. The details of these assumptions are clarified in
the following consideration. The loss probability is
the probability that the message is either marked or lost because
of overflowing the queue. We study the cases where the traffic
(offered load) $\varrho=\lambda b$ is less than, equal to and
greater than 1.

\medskip

{\bf 1.3. Advantages of the approach and methodology.}
Our approach is based on the asymptotic analysis of the loss
queueing systems in the earlier paper of the author (see Abramov [3]).
The main method is an application of modern Tauberian
theorems with remainder. For the relevant works devoted to
asymptotic analysis of the loss and controlled systems with
Poisson input see Abramov [1], [2], Tomk\'o [17] and other papers.
The asymptotic analysis of the $GI/M/1/n$ queueing system was
studied in [4], [9], [10].
The advantages of the approach of the present paper are the following.

\smallskip

{\it First}, our model is more general than the model from the
aforementioned papers: This paper discusses the case of
non-Markovian buffer model where a message contains a random batch
of packets, while the aforementioned papers studied a Markovian
model with fixed batch size.

\smallskip

{\it Second}, the work in [6], [7], [11], [13]
discusses a more difficult problem of consecutive losses,
remaining in a framework of the standard $M/M/1/n$ queueing
system. The present paper flexibly discusses the stationary losses
for a not standard queueing model with the random number of
waiting places. That queueing system belongs to the special class
of queueing systems with losses that is exactly defined below.

\smallskip

{\it Third}, our asymptotic analysis is much simpler than that
from the other papers; our final results and their
representation are simple and clear as well.

\smallskip

The traditional approach to asymptotic analysis, based on the
final value theorem for $z$ transform, enables us to obtain the
main term of asymptotic relation and in certain cases a remainder.
The modern Tauberian theorems enable us to obtain stronger
asymptotic relations using some additional assumptions.
These additional assumptions are
realistic for the queueing systems considered here, and
our asymptotic results are stronger than the earlier asymptotic
results obtained for the $M/GI/1/n$ queueing system with the aid
of the final value theorem for $z$ transform (see relation (4.15)
for its comparison with (4.14)). For some other results related to
asymptotic analysis of the $M/GI/1/n$ and $GI/M/1/n$ queueing
systems with the aid of the final value theorem see the
bibliography notes and references in Abramov [1].

\medskip

{\bf 1.4. What is the main result in this paper?} The paper
contains a number of theoretical results on the asymptotic
behavior of characteristics of the busy period of the system
(Section 4) and loss probability (Section 5). Those theoretical
results are then used to conclude the effect of adding redundant
packets in order to decrease the loss probability.

Although the theoretical results of the paper, related to the
cases where the offered load $\varrho<1$, are standard, the
conclusion about adding redundancy is extremely simple and
interesting nevertheless. Namely, the stationary loss probability
is expressed only via the probability that there is a corrupted
packet in the message. This enables us to conclude that adding a
number of redundant packets can decrease the loss probability {\it
with the rate of geometric progression} while $\varrho<1$.

Then, the case when $\varrho$ is close to 1 is very important for
the performance analysis. For example, it can be a result of the
adding a number of redundant packets when initially $\varrho<1$.
That is, to achieve a maximum decrease in the loss probability we
allow an increase in the offered load up to the critical value.

Therefore, the results on redundancy, related to the case where
$\varrho=1+\varepsilon$ ($\varepsilon>0$) is slightly greater than
1, are extremely important. The usefulness of the case
$\varrho=1+\varepsilon$ is that, it enables us to obtain more
exact conclusions on redundancy based on asymptotic results with
remainder.
Then, the usefulness of the purely theoretical case $\varrho>1$ is
that it is an intermediate result helping us to study the
transient behavior, related to the case $\varrho=1+\varepsilon$
for small $\varepsilon>0$.

\medskip

{\bf 1.5. Conclusion on adding redundant packets.} The results of
the paper enable us to make conclusions on the effect of adding
redundant packets as follows. Let $\varrho$ denote the offered
load of the system before adding a redundant packet, and let
${\mathaccent "015 \varrho}>\varrho$ be the value of offered load
after adding a redundant packet. While ${\mathaccent "015
\varrho}$ remains not greater than 1, then adding redundant
packets is profitable. It decreases the loss probability with the
rate of a geometric progression. Adding a redundant packet remains
profitable if the value ${\mathaccent "015
\varrho}=1+\varepsilon$, where $\varepsilon$ is small value of a
higher order than $p$. In some cases adding a redundant packet
decreases the loss probability even when the value $\varepsilon$
has the same order as $p$. These cases are studied in Section 6.

\medskip

{\bf 1.6. The organization of the paper.} The paper is organized
as follows. There are six sections, the first is an introduction.
In Section 2 we introduce the class of queueing systems with a
random number of waiting places and study the characteristics of
the system busy period. The results on the expectations of random
variables of the busy period (the number of processed messages,
the number of refused messages etc) are given by Lemma 2.1. In
Section 3 we present a number of auxiliary results and the
Tauberian theorems with remainder. These results are then used to
prove a number of theorems on asymptotic behavior of the
characteristics of the system given on a busy period which in turn
are given in Section 4. Section 5 presents the results on
asymptotic behavior of the loss probabilities under different
assumptions. In Section 6 we discuss adding redundancy. The
central question here is {\it how the loss probability is
decreased or increased if we add redundant packets into the
message}?

\section{Characteristics of the system given on a busy period}
\indent

The aim of this section is to to deduce the explicit
representations for characteristics of the system during a busy
period such as expected busy period, expected number of served and
lost customers during a busy period and so on. The queueing system
described in Section 1.2 is not standard, and the explicit
representation for its characteristics can not be obtained
traditionally. Therefore, below we introduce a special class of
queueing systems $\Sigma$ containing the system studied in the
paper and described in Section 1.2. It will be shown in this
section that the above characteristics are the same for all
queueing systems of the class $\Sigma$. Hence, one can take any
queueing system, a representative of class $\Sigma$, having a more
simple structure than the original system, and study it instead of
original system.

\smallskip

For the sake of convenience, we denote by ${\cal S}_1$ the system
described in Section 1.2. Let $B(x)$ be the probability
distribution function of a processing time (in the queueing
terminology - a service time), and let $\lambda$ be the parameter
of Poisson input. We also set $\varrho_j=\lambda^j\int_0^\infty
x^j\mbox{d}B(x)$, $j$=1,2,..., and $\varrho_1=\varrho$.

\medskip

In order to study the characteristics of the system ${\cal S}_1$
we introduce a set of systems $\Sigma$  containing ${\cal S}_1$ as
an element. The set $\Sigma$ is a set of $M/GI/1$ queueing systems
where $\lambda$ is the rate of Poisson input, $B(x)$ is the
probability distribution function of a service time and the family
of sequences $\{\zeta_i\}$ is more general than in ${\cal S}_1$.
Each sequence $\zeta_1$, $\zeta_2$, ... is a family of identically
distributed random variables, governing the rejection process and
having the same distribution as the random variable $\zeta$. If
this sequence is as defined in Section 1.2 then we have a
description of our system ${\cal S}_1$. In order to define the set
$\Sigma$ more exactly, we use the notation for the queueing system
${\cal S}_1$ and introduce also the following.

Let $\xi_i$ denote the number of messages in the system ${\cal
S}_1$ immediately before arrival of the $i$th message, $\xi_1=0$,
and let $s_i$ denote the number of service completions between the
$i$th and $i+1$st arrivals. It is clear that
$$
\xi_{i+1}=\xi_{i}-s_{i}+{\bf I}\{\xi_{i}\le\zeta_{i}\}, \leqno
(2.1)
$$
where the term ${\bf I}\{\xi_{i}\le\zeta_{i}\}$ indicates that
the $i$th
message is accepted, and obviously $s_i$ is not greater than
$\xi_{i}+{\bf I}\{\xi_{i}\le\zeta_{i}\}$.

Consider a new queueing system ${\cal S}$ as
above with the Poisson input rate $\lambda$, probability
distribution function of a service $B(x)$ but with the sequence
$\widetilde\zeta_1$, $\widetilde\zeta_2,...$.
Here we assume that the sequence
$\{\widetilde\zeta_i\}$ is an {\it arbitrary dependent} sequence
of random variables consisting of identically distributed
random variables as the random variable $\zeta$. Let
$\widetilde\xi_i$ denote the number of messages immediately before
arrival of the $i$th message ($\widetilde\xi_1=0$), and let
$\widetilde s_i$ denote the number of service completions between
the $i$th and $i+1$st arrivals. Thus, we assume that the initial
conditions of both queueing systems ${\cal S}_1$ and ${\cal S}$
are the same: $\xi_1=\widetilde\xi_1$.

Analogously to (2.1) we have
$$
\widetilde\xi_{i+1}=\widetilde\xi_{i}-\widetilde s_{i}+{\bf
I}\{\widetilde\xi_{i}\le\widetilde\zeta_{i}\}. \leqno (2.2)
$$

{\bf Definition.} We say that the queueing system ${\cal S}$
belongs to the set $\Sigma$ of queueing systems if ${\bf
E}\widetilde\xi_{i}={\bf E}\xi_{i}$, ${\bf E}\widetilde s_i={\bf
E}s_i$ and ${\bf
P}\{\widetilde\xi_{i}\le\widetilde\zeta_{i}\}={\bf
P}\{\xi_{i}\le\zeta_{i}\}$ for all $i\ge 1$.

\medskip

Consider an example of queueing systems belonging to the set
$\Sigma$, where the sequence $\{\widetilde\zeta_i\}$ is strictly
stationary but not ergodic.
The
example is a queueing system with
$\widetilde\zeta_1=\widetilde\zeta_2=...$, which we donate by
${\cal S}_2$.
This example below is artificial
rather than realistic, however, its main goal is to help us to
show the existence of necessary stationary queue-length
probabilities for the queueing system ${\cal S}_1$ and to obtain
the explicit representations for those probabilities as well.

For ${\cal S}_2$ we find by induction
for all $i\ge 1$ that:
$$
{\bf E}\widetilde s_{i}={\bf E}s_{i}\leqno(2.3)
$$
$$
{\bf P}\{\widetilde\xi_{i}\le\widetilde\zeta_{i}\}={\bf
P}\{\xi_{i}\le\zeta_{i}\}\leqno (2.4)
$$
and
$$
{\bf E}\widetilde\xi_{i}-{\bf E}\widetilde s_{i}+{\bf
P}\{\widetilde\xi_{i}\le\widetilde\zeta_{i}\}= {\bf E}\xi_{i}-{\bf
E}s_{i}+{\bf P}\{\xi_{i}\le\zeta_{i}\}.\leqno(2.5)
$$
Relations (2.3) - (2.5) show that the queueing system ${\cal
S}_2\in\Sigma$.

It follows from the definition that if the stationary loss
probability exists for at least one of queueing systems ${\cal
S}\in\Sigma$, then it exists for all queueing systems of $\Sigma$
and it is the same. Then, the properties of the queueing system
${\cal S}_2$ enable us to conclude similar properties of all
queueing systems belonging to the set $\Sigma$ including ${\cal
S}_1$. For example, it is not difficult to show that the expected
busy period is the same for all queueing systems of the class
$\Sigma$. Indeed, let $\widetilde A$, $\widetilde S$ and
$\widetilde R$ denote the number of arrived, served and refused
customers (because of overflowing the buffer) during a busy cycle
$\widetilde C$ respectively. We have the equations
$$
{\bf E}\widetilde A={\bf E}\widetilde S+{\bf E}\widetilde
R=\lambda {\bf E}\widetilde C,\leqno (2.6)
$$
$$
b{\bf E}\widetilde S={\bf E}\widetilde
C-\frac{1}{\lambda},\leqno(2.7)
$$
where $b$ is the expected service time. Since the loss probability
is the same for all queueing systems ${\cal S}\in\Sigma$, then the
fraction ${\bf E}\widetilde R/{\bf E}\widetilde C$ is the same for
all ${\cal S}\in\Sigma$ as well. Therefore, it follows from
equations (2.6) and (2.7) that the expected duration of a busy
period, ${\bf E}\widetilde T={\bf E}\widetilde C-\lambda^{-1}$, is
the same for all queueing systems ${\cal S}\in\Sigma$.

Recall that for that queueing system ${\cal S}_2$ we have
$\widetilde\zeta_1$=$\widetilde\zeta_2$=..., i.e. the random
variable $\zeta$ is modelled once at the initial time moment. Let
$\widetilde T_\zeta$ denote a busy period of this system. Then,
the total expectation formula enables us to write
$$
{\bf E}\widetilde
T_\zeta=\sum_{K=\zeta^{lower}}^{\zeta^{upper}}{\bf E}T_K{\bf
P}\{\zeta=K\}, \leqno (2.8)
$$
where ${\bf E}T_K$ is the expected busy period of an $M/GI/1/K$
queueing system with the same sequence of interarrival and service
times, and ${\bf P}\{\zeta=K\}={\bf P}\{\zeta_j=K\}$. In turn, the
expectation ${\bf E}T_K$ is determined from the following
recurrence relation
$$
{\bf E}T_K=\sum_{j=0}^K\pi_j{\bf E}T_{K-j+1},~~~{\bf E}T_0=b,
\leqno (2.9)
$$
$$
\pi_i=\int_0^\infty\mbox{e}^{-\lambda x}\frac{(\lambda x)^i}{i!}\mbox{d}B(x)
$$
(see Tomk\'o [17], Cooper and Tilt [12] as well as Abramov [1],
[3]) where $b$ is the expectation of a service time.

Now, let $T_\zeta$ denote a busy period for the queueing system
${\cal S}_1$. According to conclusion above that ${\bf
E}T_\zeta={\bf E}\widetilde T_\zeta$ and in view of (2.8) we have
$$
{\bf E} T_\zeta=\sum_{K=\zeta^{lower}}^{\zeta^{upper}}{\bf
E}T_K{\bf P}\{\zeta=K\}, \leqno (2.10)
$$
where ${\bf E}T_K$ are determined from (2.9).

Along with the notation $T_\zeta$ for the busy period of the
system ${\cal S}_1$, we let $I_\zeta$ be
an idle period and let $P_\zeta$, $M_\zeta$, $R_\zeta$
the characteristics of the system on a busy period:
the number of processed messages, the number of marked messages,
the number of refused messages, respectively. Here and later we use the
following terminology. The term {\it refused} message is used for
the case of overflowing the buffer. Then the term {\it lost}
message is used for the case where a message is either refused or
marked. The number of lost messages during a busy period is
denoted by $L_\zeta$. Analogously, under loss probability we mean
the probability when an arrival message is lost.

\medskip

{\bf Lemma 2.1.} {\it For the expectations ${\bf E}T_\zeta$, ${\bf
E}P_\zeta$, ${\bf E}M_\zeta$, ${\bf E}R_\zeta$ we have the
following representations:}
$$
{\bf E}P_\zeta=\frac{\lambda}{\varrho}{\bf E}T_\zeta, \leqno
(2.11)
$$
$$
{\bf E}M_\zeta=p{\bf E}P_\zeta, \leqno (2.12)
$$
$$
{\bf E}R_\zeta=(\varrho-1){\bf E}P_\zeta+1. \leqno(2.13)
$$

{\bf Proof}. Relations (2.11) and (2.12) follow immediately from
the Wald's identity. In order to prove (2.13) note that the number
of arrivals during a busy cycle equals to the number of processed
messages during a busy period plus the number of refused messages
during a busy period (see relation (2.6)). According to the Wald's
identity the expected number of arrivals during a busy cycle
equals $\lambda({\bf E}T_\zeta+{\bf E}I_\zeta)$. Therefore taking
into account that ${\bf E}I_\zeta=\lambda^{-1}$ from (2.11) we
have
$$
{\bf E}R_\zeta=(\varrho-1){\bf E}P_\zeta+1,
$$
and the result is proved. \\

For the alternative proof of (2.13) see Abramov [3]. (See also the
proof in [5].)

\section{Auxiliary results. Tauberian theorems with remainder}
\indent

It is seen from relations (2.10) and (2.9) and Lemma 2.1 that the
characteristics of the system during a busy period can be
studied in a framework of the recurrence relation
$$
Q_k=\sum_{i=0}^k r_iQ_{k-i+1}, \leqno(3.1)
$$
where $r_i$ are nonnegative numbers, $r_0+r_1+...=1$, $r_0>0$, and $Q_0>0$ is an
arbitrary real number.
Below we recall a number of results
on asymptotic behavior of that sequence (3.1).

The known results on representation (3.1) are asymptotic theorems
by Tak\'acs [16]. Lemma 3.1 below joins two results by Tak\'acs
[16]: Theorem 5 on p. 22 and relation (35) on p. 23. The results
of Tak\'acs [16] were then developed by Postnikov [14, Section 25]
(see Lemma 3.2 and Lemma 3.3 below).

Let $r(z)=\sum_{i=0}^\infty r_iz^i$, $|z|\le 1$,
$\gamma_m=r^{(m)}(1-0)=\lim_{z\uparrow 1}r^{(m)}(z)$
($r^{(m)}(z)$ is the $m$th derivative of $r(z)$). Note that if denote
$Q(z)=\sum_{i=0}^\infty Q_iz^i$
then it follows from (3.1) that
$$Q(z)=\frac{Q_0r(z)}{r(z)-z}.$$

{\bf Lemma 3.1.} (Tak\'acs [16]). {\it If $\gamma_1<1$ then
$$
\lim_{k\to\infty}Q_k=\frac{Q_0}{1-\gamma_1}. \leqno (3.2)
$$
If $\gamma_1=1$ and $\gamma_2<\infty$ then
$$
\lim_{k\to\infty}\frac{Q_k}{k}=\frac{2Q_0}{\gamma_2}. \leqno (3.3)
$$
If $\gamma_1>1$ then
$$
\lim_{k\to\infty}\Big(Q_k-\frac{Q_0}{\delta^k[1-r'(\delta)]}\Big)=\frac
{Q_0}{1-\gamma_1}, \leqno (3.4)
$$
where $\delta$ is the least (absolute) root of equation} $z=r(z)$.\\

{\bf Lemma 3.2.} (Postnikov [14]). {\it Let $\gamma_1=1$,
$\gamma_3<\infty$. Then as} $k\to\infty$
$$
Q_k=\frac{2Q_0}{\gamma_2}k+O(\log k). \leqno (3.5)
$$

{\bf Lemma 3.3.} (Postnikov [14]). {\it Let $\gamma_1=1$,
$\gamma_2<\infty$ and $r_0+r_1<1$. Then as} $k\to\infty$
$$
Q_{k+1}-Q_k=\frac{2Q_0}{\gamma_2}+o(1). \leqno (3.6)
$$

\section{Asymptotic results for characteristics of the system
during a busy period}
\indent

This section provides a number of results on asymptotic behavior
of characteristics of the system. The first three theorems are
related to the case as $N$ increases to infinity, where the cases
$\varrho<1$, $\varrho=1$ and $\varrho>1$ are considered. The next
two theorems discuss the case of the value $\varrho$ is close to
the critical value 1, and as $N\to\infty$, it tends to 1. The last
theorem of this section, Theorem 4.6, provides the asymptotic
result for the special case when the number of packets in a
message is a constant value.

\smallskip

Let us now study the asymptotic behavior of the expectations ${\bf
E}P_\zeta$, ${\bf E}M_\zeta$ and ${\bf E}R_\zeta$. We write
$\zeta=\zeta(N)$ pointing out the dependence on parameter $N$. As
the buffer size $N$ increases to infinity then both
$\zeta^{lower}$ and $\zeta^{upper}$ tend to infinity, and together
with them $\zeta(N)$ a.s. tends to infinity. Then we have the
following.

\medskip

{\bf Theorem 4.1.} {\it If $\varrho<1$ then
$$
\lim_{N\to\infty}{\bf E}P_{\zeta(N)}=\frac{1}{1-\varrho}. \leqno
(4.1)
$$
If $\varrho=1$ and $\varrho_2<\infty$ then
$$
\lim_{N\to\infty}\frac{{\bf E}P_{\zeta(N)}}{{\bf
E}\zeta(N)}=\frac{2}{\varrho_2}. \leqno (4.2)
$$
If $\varrho>1$ then
$$
\lim_{N\to\infty}\Big[{\bf E}P_{\zeta(N)}-\frac{1}{{\bf
E}\varphi^{\zeta(N)}[1+\lambda\beta'(\lambda-\lambda\varphi)]}\Big]=\frac{1}{1-\varrho},\leqno
(4.3)
$$
where $\beta(z)=\int_0^\infty\mbox{e}^{-zx}\mbox{d}B(x)$ and
$\varphi$ is the least (absolute) root of functional equation}
$z-\beta(\lambda-\lambda z)
=0$.\\

{\bf Proof}. From (2.9), (2.10) and (2.11) we have
$$
{\bf E}P_\zeta=\sum_{K=\zeta^{lower}}^{\zeta^{upper}}{\bf
E}P_K{\bf P}\{\zeta=K\},
$$
where
$$
{\bf E}P_K=\sum_{j=0}^K\pi_j{\bf E}P_{K-j+1},~~~{\bf E}P_0=1,
$$
$$
\pi_j=\int_0^\infty\mbox{e}^{-\lambda x}\frac{(\lambda x)^j}{j!}\mbox{d}B(x).
$$
Then applying Lemma 3.1 we have the following. In the case
$\varrho<1$, taking into account that $\zeta(N){\buildrel
a.s.\over\to}\infty$ as $N\to\infty$, we obtain
$$
\lim_{N\to\infty}{\bf E}P_{\zeta(N)}=\lim_{N\to\infty}{\bf
E}P_N=\frac{1}{1-\varrho}.
$$
Relation (4.1) is proved.

In the case $\varrho_2<\infty$ and $\varrho=1$ we have
$$
\lim_{N\to\infty}\frac{{\bf
E}P_{\zeta(N)}}{N}=\lim_{N\to\infty}\frac{1}{N}\sum_{K=\zeta^{lower}}^{\zeta^{upper}}{\bf
P} \{\zeta(N)=K\}{\bf E}P_K
$$
$$
=\lim_{N\to\infty}\frac{1}{N}\sum_{K=\zeta^{lower}}^{\zeta^{upper}}
K{\bf P}\{\zeta(N)=K\}\frac{2}{\varrho_2}
=\frac{2}{\varrho_2}\lim_{N\to\infty}\frac{{\bf E}\zeta(N)}{N}.
$$
Therefore,
$$
\lim_{N\to\infty}\frac{{\bf E}P_{\zeta(N)}}{{\bf
E}\zeta(N)}=\frac{2}{\varrho_2},
$$
and relation (4.2) is proved.

In the case $\varrho>1$ for large $N$ we obtain
$$
{\bf E}P_{\zeta(N)}=\sum_{K=\zeta^{lower}}^{\zeta^{upper}} {\bf
P}\{\zeta(N)=K\}{\bf E}P_K
$$
$$
=\sum_{K=\zeta^{lower}}^{\zeta^{upper}} {\bf
P}\{\zeta(N)=K\}\frac{1}{\varphi^{K}[1+\lambda
\beta'(\lambda-\lambda\varphi)]}+\frac{1}{1-\varrho}+o(1)
$$
$$
=\frac{1}{{\bf
E}\varphi^{\zeta(N)}[1+\lambda\beta'(\lambda-\lambda\varphi)]}
+\frac{1}{1-\varrho}+o(1).
$$
Therefore,
$$
\lim_{N\to\infty}\Big[{\bf E}P_{\zeta(N)}-\frac{1}{{\bf
E}\varphi^{\zeta(N)}[1+\lambda\beta'(\lambda-\lambda\varphi)]}\Big]=\frac{1}{1-\varrho},
$$
and relation (4.3) is proved. Theorem 4.1 is completely proved.\\

{\bf Theorem 4.2.} {\it If $\varrho=1$ and $\varrho_3<\infty$
then}
$$
{\bf E}P_{\zeta(N)}=\frac{2}{\varrho_2}{\bf E}\zeta(N)+O(\log N).
\leqno (4.4)
$$

{\bf Proof}.
Applying Lemma 3.2, for large $N$ we have
$$
{\bf E}P_{\zeta(N)}= \sum_{K=\zeta^{lower}}^{\zeta^{upper}}{\bf P}
\{\zeta(N)=K\}{\bf E}P_K =\sum_{K=\zeta^{lower}}^{\zeta^{upper}}
K{\bf P} \{\zeta(N)=K\}\frac{2}{\varrho_2}+O\{{\bf E}[\log
\zeta(N)]\}
$$
$$
=\frac{2}{\varrho_2}{\bf E}\zeta(N)+O\{{\bf E}[\log \zeta(N)]\}
$$
$$
=\frac{2}{\varrho_2}{\bf E}\zeta(N)+O(\log N).~~~~~~~
$$
and we obtain relation (4.4). Theorem is proved.

\medskip

In turn for ${\bf E}R_{\zeta}$ we have

\medskip

{\bf Theorem 4.3.} {\it If $\varrho<1$ then
$$
\lim_{N\to\infty}{\bf E}R_{\zeta(N)}=0. \leqno (4.5)
$$
If $\varrho=1$ then for all $N\ge 0$
$$
{\bf E}R_{\zeta(N)}=1. \leqno (4.6)
$$
If $\varrho>1$ then}
$$
\lim_{N\to\infty}\Big[{\bf E}R_{\zeta(N)}-\frac{\varrho-1}{{\bf
E}\varphi^{\zeta(N)}[1+\lambda\beta'(\lambda-\lambda\varphi)]}\Big]=0,\leqno
(4.7)
$$

\medskip

{\bf Proof}. The proof of this theorem is analogous to that of the
proof of Theorem 4.1. It follows by application of Lemma 3.1 and
relation (2.13) of Lemma 2.1.

\medskip

{\bf Theorem 4.4.} {\it Let $\varrho=1+\varepsilon$,
$\varepsilon>0$, and $\varepsilon \zeta(N)\to C>0$ a.s. as
$\varepsilon\to 0$ and $N\to\infty$. Assume also that
$\varrho_3=\varrho_3(N)$ is a bounded sequence, and there exists
${\mathaccent "7E \varrho}_2=\lim_{N\to\infty}\varrho_2(N)$. Then}
$$
{\bf E}P_{\zeta(N)}=\frac{\mbox{e}^{2C/{\mathaccent "7E
\varrho}_2}-1}{\varepsilon}+O(1), \leqno (4.8)
$$
$$
{\bf E}R_{\zeta(N)}=\mbox{e}^{2C/{\mathaccent "7E
\varrho}_2}+o(1). \leqno (4.9)
$$

\medskip

{\bf Proof}. It was shown in Subhankulov [15, p. 326], that if
$\varrho=1+\varepsilon$, $\varepsilon>0$, $\varepsilon\to 0$,
$\varrho_3(N)$ is a bounded sequence, and there exists
${\mathaccent "7E \varrho}_2=\lim_{N\to\infty}\varrho_2(N)$, then
$$
\varphi=1-\frac{2\varepsilon}{{\mathaccent "7E
\varrho}_2}+O(\varepsilon^2). \leqno (4.10)
$$
Applying (4.10) after some algebra we have
$$
1+\lambda\beta'(\lambda-\lambda\varphi)=\varepsilon+O(\varepsilon^2).
\leqno (4.11)
$$
Then the statements of the theorem follow by applying expansions
(4.10) and (4.11) to (4.3) and (4.7).

\medskip

{\bf Theorem 4.5.} {\it Let $\varrho=1+\varepsilon$,
$\varepsilon>0$, and $\varepsilon \zeta(N)\to 0$ as
$\varepsilon\to 0$ and $N\to\infty$. Assume also that
$\varrho_3=\varrho_3(N)$ is a bounded sequence, and there exists
${\mathaccent "7E \varrho}_2=\lim_{N\to\infty}\varrho_2(N)$. Then}
$$
{\bf E}P_{\zeta(N)}=\frac{2}{\widetilde\varrho_2}{\bf
E}\zeta(N)+O(1), \leqno (4.12)
$$
$$
{\bf E}R_{\zeta(N)}=1+o(1). \leqno (4.13)
$$

{\bf Proof.} The results follow by expanding (4.8) and (4.9) for
small $C$.

\medskip

{\bf Special case.} If each message contains the same number of
packets, say $l$, then we have the usual $M/GI/1/n$ queueing
system, where $n=\lfloor N/l \rfloor$ is the integer part of
$N/l$. For that queueing system all the results in
Theorems 4.1-4.5 hold, by replacing $\zeta(N)$ (or ${\bf
E}\zeta(N)$) by $n$.

For example, asymptotic relation (4.7) looks
$$
\lim_{n\to\infty}\Big({\bf
E}R_n-\frac{\varrho-1}{\varphi^n[1+\lambda\beta'(\lambda-\lambda\varphi)]}\Big)=0.
\leqno (4.14)
$$
Notice that using the final value theorem for $z$ transform,
Azlarov and Tahirov [8] obtained the estimation
$$
{\bf
E}R_n=\frac{\varrho-1}{\varphi^n[1+\lambda\beta'(\lambda-\lambda\varphi)]}
\Big[1+O\Big(\frac{2\varphi}{1+\varphi}\Big)^n\Big],
\leqno (4.15)
$$
weaker than (4.14).

The theorem below is related to the case of the usual queueing
systems only, when the number of packets in a message is fixed.
Namely, we have

\medskip

{\bf Theorem 4.6.} {\it If $\varrho=1$ and $\varrho_2<\infty$ then
$$
{\bf E}P_{n+1}-{\bf E}P_n=\frac{2}{\varrho_2}+o(1),~~~n\to\infty,
\leqno (4.16)
$$
where the index $n+1$ says that $P_{n+1}$ is the number of
processed
messages during a busy period of the $M/GI/1/n+1$ queueing system.}\\

{\bf Proof}. The result will follow from Lemma 3.3 if we show that
$\beta(\lambda)-\lambda\beta'(\lambda)<1$. Taking into account
that for each $\lambda>0$
$$
\sum_{i=0}^\infty \frac{(-\lambda)^i}{i!}\beta^{(i)}(\lambda)=
\sum_{i=0}^\infty \int_0^\infty\mbox{e}^{-\lambda x}\frac{(\lambda
x)^i}{i!} \mbox{d}B(x)=
\int_0^\infty\sum_{i=0}^\infty\mbox{e}^{-\lambda x}\frac{(\lambda
x)^i}{i!} \mbox{d}B(x)=1, \leqno (4.17)
$$
and all terms
$$
\pi_i=\frac{(-\lambda)^i}{i!}\beta^{(i)}(\lambda)
$$
are nonnegative, from (4.17) we find that
$$
\beta(\lambda)-\lambda\beta'(\lambda)\le 1. \leqno (4.18)
$$
Thus, the required statement will be proved if we show that for
some $\lambda_0>0$ the equality
$$
\beta(\lambda_0)-\lambda_0\beta'(\lambda_0)=1 \leqno (4.19)
$$
is not a case. Indeed, since the function
$\beta(\lambda)-\lambda\beta'(\lambda)$ is an analytic function
then according to the maximum absolute value principle for
analytic function $\beta(\lambda)-\lambda\beta'(\lambda)=1$ holds
for all $\lambda>0$. Therefore identity (4.19) means that
$\pi_i=0$ for all $i\ge 2$ and for all $\lambda>0$. Therefore,
(4.19) is valid if and only if $\beta(\lambda)$ is a linear
function, i.e. $\beta(\lambda)=c_0+c_1\lambda$, $c_0$ and $c_1$
are some constants. However, since $|\beta(\lambda)|\le 1$ we
obtain $c_0=1$ and $c_1=0$, and $\beta(\lambda)\equiv 1$. This is
the trivial case where the probability distribution function
$B(x)$ is concentrated in point 0. Therefore (4.19) is not a case,
and $\beta(\lambda)-\lambda\beta'(\lambda)<1$. The theorem is
proved.

\section{Asymptotic theorems for the loss probabilities}
\indent

In this section we study the asymptotic behavior of the loss
probability by using renewal arguments. The results of this
section correspond to those of the previous section. We discuss
the behavior of the system for the same cases as $N\to\infty$, as
well when the parameter $\rho$ is close to the critical value 1
and tends to 1 as $N\to\infty$. The theorems of this section are
important for our conclusion on adding redundancy, which is given
in the next section.

\smallskip

According to renewal arguments the loss probability is determined
as
$$\Pi_\zeta=\frac{{\bf E}L_\zeta}{{\bf E}R_\zeta+{\bf E}P_\zeta}=
\frac{{\bf E}R_\zeta+{\bf E}M_\zeta}{{\bf E}R_\zeta+{\bf
E}P_\zeta} =\frac{{\bf E}R_\zeta+p{\bf E}P_\zeta}{{\bf
E}R_\zeta+{\bf E}P_\zeta}. \leqno (5.1)
$$
(Recall that $L_\zeta$ is the number of lost messages during a busy period.)\\

{\bf Theorem 5.1.} {\it If $\varrho<1$
$$
\lim_{N\to\infty}\Pi_{\zeta(N)}=p. \leqno (5.2)
$$
(Recall that $p$ is the probability that a message is erroneous
because one of its packets is corrupted.)

Limiting relation (5.2) is also valid when $\varrho=1$ and
$\varrho_2<\infty$.

If $\varrho>1$ then}
$$
\Pi_{\zeta(N)}=\frac{p+\varrho-1}{\varrho}~\frac{(\varrho-1)+
p[1+\lambda\beta'(\lambda-\lambda\varphi)]{\bf
E}\varphi^{\zeta(N)}}
{(\varrho-1)+[1+\lambda\beta'(\lambda-\lambda\varphi)]{\bf
E}\varphi^{\zeta(N)}}+o({\bf E}\varphi^{\zeta(N)}). \leqno (5.3)
$$

{\bf Proof}.
The proof follows from Theorems 4.1 and 4.3.\\

{\bf Theorem 5.2.} {\it If $\varrho=1$ and $\varrho_3<\infty$ then
as} $N\to\infty$
$$
\Pi_{\zeta(N)}=p+\frac{(1-p)\varrho_2}{2{\bf
E}\zeta(N)}+O\Big(\frac{\log N}{N^2}\Big). \leqno (5.4)
$$

{\bf Proof}. From (5.1) we have
$$
\Pi_{\zeta(N)}=\frac{{\bf E}R_{\zeta(N)}}{{\bf E}R_{\zeta(N)}+{\bf
E}P_{\zeta(N)}}+ \frac{p{\bf E}P_{\zeta(N)}}{{\bf
E}R_{\zeta(N)}+{\bf E}P_{\zeta(N)}} \leqno (5.5)
$$
$$
=\frac{1}{1+{\bf E}P_{\zeta(N)}}+
\frac{p{\bf E}P_{\zeta(N)}}{1+{\bf E}P_{\zeta(N)}}.
$$
As $N\to\infty$ from Theorem 4.2 we obtain:
$$
\frac{1}{1+{\bf E}P_{\zeta(N)}}=\frac{\varrho_2}{2{\bf
E}\zeta(N)}+ O\Big(\frac{\log N}{N^2}\Big), \leqno (5.6)
$$
$$
\frac{p{\bf E}P_{\zeta(N)}}{1+{\bf
E}P_{\zeta(N)}}=p-\frac{p\varrho_2}{2{\bf
E}\zeta(N)}+O\Big(\frac{\log N}{N^2}\Big). \leqno (5.7)
$$
Combining these two asymptotic relations (5.6) and (5.7) we obtain
the statement of Theorem 5.2. Theorem 5.2 is proved.

\medskip

{\bf Note.} Under assumptions of Theorem 5.2 assume additionally that $p\to 0$. If $pN\to C>0$ then
$$
\Pi_{\zeta(N)}=\frac{C}{N}+\frac{\varrho_2}{2{\bf
E}\zeta(N)}+O\Big(\frac{\log N}{N^2}\Big).
$$
If $pN\to 0$ then
$$
\Pi_{\zeta(N)}=\frac{\varrho_2}{2{\bf
E}\zeta(N)}+O\Big(p+\frac{\log N}{N^2}\Big).
$$

\medskip

The theorem below also assumes that $p\to 0$. Our result here is
the following.

\medskip

{\bf Theorem 5.3.} {\it Let $\varrho=1+\varepsilon$,
$\varepsilon>0$, and $\varepsilon \zeta(N)\to C>0$ as
$\varepsilon\to 0$ and $N\to\infty$, and $p\to 0$. Assume also
that $\varrho_3=\varrho_3(N)$ is a bounded sequence, and there
exists ${\mathaccent "7E
\varrho}_2=\lim_{n\to\infty}\varrho_2(N)$.

\medskip

$(i)$  If $p/\varepsilon\to D\ge 0$ then we have
$$
\Pi_{\zeta(N)}=\Big(D+\frac{\mbox{e}^{2C/\widetilde\varrho_2}}{\mbox{e}^{2C/\widetilde\varrho_2}-1}\Big)
\varepsilon+o(\varepsilon).
\leqno (5.8)
$$

$(ii)$ If $p/\varepsilon\to\infty$ then we have}
$$
\Pi_{\zeta(N)}=p+O(\varepsilon). \leqno (5.9)
$$

{\bf Proof.} In the case $(i)$ we have
$$
p{\bf E}P_{\zeta(N)}+{\bf
E}R_{\zeta(N)}=(D+1)\mbox{e}^{2C/\widetilde\varrho_2}-D+o(1),
\leqno (5.10)
$$
and
$$
{\bf E}P_{\zeta(N)}+{\bf
E}R_{\zeta(N)}=\frac{\mbox{e}^{2C/\widetilde\varrho_2}-1}{\varepsilon}+O(1).
\leqno (5.11)
$$
Therefore from (5.10) and (5.11) we have
$$
\Pi_{\zeta(N)}=\Big(D+\frac{\mbox{e}^{2C/\widetilde\varrho_2}}{\mbox{e}^{2C/\widetilde\varrho_2}-1}\Big)\varepsilon+o(\varepsilon),
$$
and relation (5.8) is proved.

In the case $(ii)$ we have
$$
p{\bf E}P_{\zeta(N)}+{\bf
E}R_{\zeta(N)}=\frac{pc}{\varepsilon}+O(1), \leqno (5.12)
$$
and
$$
{\bf E}P_{\zeta(N)}+{\bf
E}R_{\zeta(N)}=\frac{c}{\varepsilon}+O(1), \leqno (5.13)
$$
where
$c=\exp(2C/\widetilde\varrho_2)/(\exp(2C/\widetilde\varrho_2)-1)$.
Relation (5.9) follows.

\medskip

{\bf Theorem 5.4.} {\it Let $\varrho=1+\varepsilon$,
$\varepsilon>0$, and $\varepsilon \zeta(N)\to 0$ as
$\varepsilon\to 0$ and $N\to\infty$, and $p\to 0$. Assume also
that $\varrho_3=\varrho_3(N)$ is a bounded sequence, and there
exists ${\mathaccent "7E
\varrho}_2=\lim_{n\to\infty}\varrho_2(N)$.

\medskip

$(i)$  If $p/\varepsilon\to D\ge 0$ then we have
$$
\Pi_{\zeta(N)}=p+\frac{\widetilde\varrho_2}{2{\bf E}\zeta(N)
}+o\Big(\frac{1}{N}\Big). \leqno (5.14)
$$

$(ii)$ If $p/\varepsilon\to\infty$ then we have (5.9).}

\medskip

{\bf Proof.} The proof of (5.14) follows by expanding (5.8) for
small $C$. The proof in case $(ii)$ trivially follows from (5.12)
and (5.13).

\medskip

{\bf Special case.} In the case where each message contains
exactly $l$ packets, $n=\lfloor N/l \rfloor$, we obtain the
following

\medskip

{\bf Theorem 5.5.} {\it If $\varrho=1$ and $\varrho_2<\infty$ then
as} $n\to\infty$
$$
\Pi_{n+1}-\Pi_n=\frac{\frac{1}{n(n+1)}~\frac{2}{\varrho_2}(p-1)}
{\Big(\frac{2}{\varrho_2}+\frac{1}{n+1}\Big)\Big(\frac{2}{\varrho_2}+\frac{1}{n}
\Big)}+o\Big(\frac{1}{n^2}\Big). \leqno (5.15)
$$

{\bf Proof}. The proof follows by application of Theorem 4.5 and
taking into account the fact that ${\bf E}R_n=1$ for all $n\ge 0$
(see [3] or Lemma 2.1)

\section{Adding redundant packets}
\indent

We now investigate the effect of adding redundant packets. We
assume that adding a redundant packet to the message decreases the
probability $p$ that a message is corrupted and increases the
offered load and the number of packets in a message. The new
parameters of the system after adding a redundant packet are
denoted by adding the symbol $\; \breve{ } \;$ above. For example,
${\mathaccent "015 p}$ is a probability that a message contains a
corrupted packet and ${\mathaccent "015 \varrho}$ is the offered
load. It follows from Theorem 5.1 that if ${\mathaccent "015
\varrho}\le 1$ the stationary loss probability coincides with
${\mathaccent "015 p}$. This means that if adding a redundant
packet to the message decreases the probability $p$ by $\gamma$
times, then the same effect is achieved with the loss probability.
Thus, adding a number of redundant packets while $\varrho<1$ can
decrease the loss probability geometrically.

In the case where both $\varrho>1$ and ${\mathaccent "015
\varrho}> 1$ the adding a redundant packet to the message changes
the stationary loss probability approximately to
$$
\frac{\varrho ({\mathaccent "015 p}+{\mathaccent "015
\varrho}-1)}{{\mathaccent "015 \varrho}(p+\varrho-1)}.
$$
In practice the values $p$ and ${\mathaccent "015 p}$ are
small, and even if adding redundant packets can slightly decrease
the stationary loss probability, the effect of that action is not
considerable.

The case where $\varrho<1$ and ${\mathaccent "015 \varrho}> 1$ is
especially interesting if ${\mathaccent "015 \varrho}=1+\delta$,
and $\delta$ is a small value. For example, if $\delta$ is so
small that both $\delta\zeta(N)$ and $\delta/p$ are also
negligible, then a redundant packet decreases the loss probability
approximately by the same amount as in the case when both
$\varrho<1$ and ${\mathaccent "015 \varrho}< 1$. However, if
$\delta$ is of the same order as $p$ or $1/\zeta(N)$ the special
analysis based on the corresponding cases of Theorems 5.3 and 5.4
is necessary. Here we do not provide the details.

Let us consider the cases when both $\varrho>1$ and ${\mathaccent
"015 \varrho}> 1$, where $\varrho=1+\varepsilon$ and ${\mathaccent
"015 \varrho}=1+{\mathaccent "015 \varepsilon}$, and $\varepsilon$
and ${\mathaccent "015 \varepsilon}$ are small values as in
Theorem 5.3, both satisfying ($i$). Then the stationary loss
probability is changed approximately to
$$
\frac{\mbox{e}^{2C/\widetilde\varrho_2}-1}
{\mbox{e}^{2{\mathaccent "015 C}/{\mathaccent "015 \varrho}_2}-1}~
\frac{(\mbox{e}^{2{\mathaccent "015 C}/{\mathaccent "015
\varrho}_2}-1){\mathaccent "015 p}+\mbox{e}^ {2{\mathaccent "015
C}/{\mathaccent "015 \varrho_2}}{\mathaccent "015 \varepsilon}}
{(\mbox{e}^{2C/\widetilde\varrho_2}-1)p+\mbox{e}^{2C/\widetilde\varrho_2}\varepsilon}~.
\leqno (6.1)
$$
times.

For the sake of simplicity let us assume that ${\mathaccent "015
C}/{\mathaccent "015 \varrho}_2=C/\widetilde\varrho_2$. Then (6.1)
reduces to
$$
\frac{(\mbox{e}^{2C/\widetilde\varrho_2}-1){\mathaccent "015
p}+\mbox{e}^ {2C/\widetilde\varrho_2}{\mathaccent "015
\varepsilon}}
{(\mbox{e}^{2C/\widetilde\varrho_2}-1)p+\mbox{e}^{2C/\widetilde\varrho_2}\varepsilon}.
\leqno (6.2)
$$

If we assume that
$$
p-{\mathaccent "015
p}=\frac{\mbox{e}^{2C/\widetilde\varrho_2}}{\mbox{e}^{2C/\widetilde\varrho_2}-1}
({\mathaccent "015 \varepsilon}-\varepsilon),
$$
then the stationary loss probability remains approximately at the
same value, and if
$$
p-{\mathaccent "015
p}>\frac{\mbox{e}^{2C/\widetilde\varrho_2}}{\mbox{e}^{2C/\widetilde\varrho_2}-1}
({\mathaccent "015 \varepsilon}-\varepsilon),
$$
then the stationary loss probability decreases, otherwise if
$$
p-{\mathaccent "015
p}<\frac{\mbox{e}^{2C/\widetilde\varrho_2}}{\mbox{e}^{2C/\widetilde\varrho_2}-1}
({\mathaccent "015 \varepsilon}-\varepsilon),
$$
then the stationary loss probability increases.

\section*{Acknowledgement}
\indent

The author thanks Professor Moshe Sidi (Technion) for sending him
the files of related papers. The author thanks also the anonymous
referees and associate editor for a number of valuable comments.

\section*{References}
\indent

[1] {\sc V.M.Abramov}, {\em Investigation of a Queueing System
with Service Depending on Queue Length}. Donish, Dushanbe, 1991.
(In Russian).

\smallskip
[2] {\sc V.M.Abramov}, {\em Asymptotic theorems for one queueing
system with refusals}. Kibernetika (2) (1991) pp. 123-124. (In
Russian).

\smallskip
[3] {\sc V.M.Abramov}, {\em On a property of a refusals stream}.
J. Appl. Probab. 34 (1997) pp. 800-805.

\smallskip
[4] {\sc V.M.Abramov}, {\em Asymptotic analysis of the $GI/M/1/n$
loss system as $n$ increases to infinity}. Ann. Operat. Res. 112
(2002) pp. 35-41.

\smallskip
[5] {\sc V.M.Abramov}, {\em On losses in $M^X/GI/1/n$ queues}. J.
Appl. Probab. 38 (2001) pp. 1079-1080.

\smallskip
[6] {\sc O.Ait-Hellal, E.Altman, A.Jean-Marie and I.A.Kurkova},
{\em On loss probabilities in presence of redundant packets and
several traffic sources}. Perform. Eval. 36-37 (1999) pp. 485-518.

\smallskip
[7] {\sc E.Altman and A.Jean-Marie}, {\em Loss probabilities for
messages with redundant packets feeding a finite buffer}. IEEE J.
Select. Areas Commun. 16 (1998) pp. 778-787.

\smallskip
[8] {\sc T.A.Azlarov and A.Tahirov}, {\em Limit theorems for
single-server queueing system with finite number of waiting
places}. Proc. USSR Acad. Sci. Engineering Cybern. (5) (1974) pp.
53-57.


\smallskip
[9] {\sc B.D.Choi and B.Kim}, {\em Sharp results on convergence
rates for the distribution of $GI/M/1/K$ queues as $K$ tends to
infinity}. J. Appl. Probab. 37 (2000) pp. 1010-1019.

\smallskip
[10] {\sc B.D.Choi, B.Kim and I.-S.Wee}, {\em Asymptotic behavior
of loss probability in the $GI/M/1/K$ queue as $K$ tends to
infinity}. Queueing Systems 36 (2000) pp. 437-442.

\smallskip
[11] {\sc I.Cidon, A.Khamisy and M.Sidi},  {\em Analysis of packet
loss processes in high-speed networks}. IEEE Trans. Inform. Theory
39 (1993) pp. 98-108.

\smallskip
[12] {\sc R.B.Cooper and B.Tilt}, {\em On the relationship between
the distribution of the maximum queue-length in the $M/G/1$ queue
and the mean busy period in the $M/G/1/n$ queue}. J. Appl. Probab.
13 (1976) pp. 195-199.

\smallskip
[13] {\sc O.Gurewitz, M.Sidi and I.Cidon}, {\em The ballot theorem
strikes again: Packet loss process distribution}. IEEE Trans.
Inform. Theory 46 (2000) pp. 2588-2595.

\smallskip
[14] {\sc A.G.Postnikov}, {\em Tauberian Theory and Its
Application}. Proc. Steklov Math. Inst. 144, No.2 (1979), pp.
1-148, and Amer. Math. Soc. (1980), pp. 1-138.

\smallskip
[15] {\sc M.A.Subhankulov}, {\em Tauberian Theorems with
Remainder}. Nauka, Moscow, 1976. (In Russian).

\smallskip
[16] {\sc L.Tak\'acs}, {\em Combinatorial Methods in the Theory of
Stochastic Processes}. John Wiley, New York, 1967.

\smallskip
[17] {\sc J.Tomk\'o}, {\em One limit theorem in the queueing
problem as input rate increases infinitely}. Studia Sci. Math.
Hungarica 2 (1967) pp. 447-454. (In Russian.)

\medskip

\end{document}